\documentclass[12pt]{amsart}
\usepackage{amsmath,amssymb,amsfonts,amsthm,amsopn}
\usepackage{graphics}

\newcommand{\ft}{Fourier transform}

\newtheorem{tm}{Theorem}[section]
\newtheorem{lemma}[tm]{Lemma}
\newtheorem{prop}[tm]{Proposition}

\newtheorem{theorem}{Theorem}[section]

\newtheorem{proposition}[theorem]{Proposition}
\newtheorem{remark}[theorem]{Remark}

\newcommand{\beqa}{\begin{eqnarray*}}
\newcommand{\eeqa}{\end{eqnarray*}}

\newcommand{\field}[1]{\mathbb{#1}}
\newcommand{\bR}{\field{R}}        
 \def\cF{\mathcal{F}}              
 \def\cS{\mathcal{S}}
 \def\cD{\mathcal{D}}

 \def\cC{\mathcal{C}}

\def\a{\aleph}

\def\rd{\bR^d}

\def\lrd{L^2(\rd)}

\def\intrd{\int_{\rd}}

\def\R{\right)}

\def\<{\left<}
\def\>{\right>}

\def\inv{^{-1}}

\def\mv1{M_v^1}

\def\noi{\noindent}

\def\o{\omega}
\def\a{\alpha}

\def\ZZ{\mathbb{Z}}

\def\R{\mathbb{R}}
\def\Ren{\mathbb{R}^d}
\def\Renn{\mathbb{R}^{2d}}

\def\sch{\mathcal{S}}

\def\Fur{\mathcal{F}}

\def\f{\varphi}

\def\Sn2{S_{2}(L^{2}(\Ren))}
\def\S1{S_{1}(L^{2}(\Ren))}
\def\sig00{\sigma_{0,0}}

\begin{document}

\title[Strichartz Estimates for the Schr\"odinger equation]{Strichartz Estimates in Wiener Amalgam Spaces for the Schr\"odinger equation}

\author{Elena Cordero and Fabio Nicola}
\address{Department of Mathematics,  University of Torino, Italy}
\address{Dipartimento di Matematica, Politecnico di Torino, Italy}



\date{}

\begin{abstract} We study the dispersive properties of the
 Schr\"o\-din\-ger equation. Precisely, we look for estimates which give a control of the local regularity and decay at infinity {\it separately}.
The Banach spaces that allow such a treatment are the Wiener amalgam spaces, and  Strichartz-type estimates are proved in
this framework. These estimates improve some of the
classical ones in the case of large time.
\end{abstract}

\maketitle

\section{Introduction}
The study of space-time integrability properties of the solution of the Cauchy problem for the Schr\"odinger equation
\begin{equation}\label{cp}
\begin{cases}
i\partial_t u+\Delta u=0,\\
u(0,x)=u_0(x),
\end{cases}
\end{equation}
with $(t,x)\in\bR\times\bR^d$, $d\geq1$, has been pursued by many authors in the last thirty years. The celebrated
homogeneous Strichartz estimates \cite{GinibreVelo92,Kato70,keel,Yajima87} for the solution
$u(t,x)=\left(e^{it\Delta}u_0\right)(x)$ read
\begin{equation}\label{S1}
\|e^{it\Delta} u_0\|_{L^q_t L^r_x}\lesssim \|u_0\|_{L^2_x},
\end{equation}
for $q\geq2$, $r\geq 2$, with
$2/q+d/r=d/2$, $(q,r,d)\not=(2,\infty,2)$, i.e., for $(q,r)$ {\it Schr\"odinger admissible}. Here, as usual, we set
\[
\|F\|_{L^q_tL^r_x}=\left(\int\|F(t,\cdot)\|^q_{L^r_x}\,dt\right)^{1/q}.
\]
As a matter of fact, these estimates express  a gain of local $x$-regularity of the solution $u(t,\cdot)$, and a decay of its $L^r_x$-norm,
 both in some $L^q_t$-averaged sense.
\par
In this paper we study similar estimates in spaces which, unlike
the $L^p$ spaces, control  the local and global behaviour of a
function independently (for example, spaces whose functions are
locally in some $L^p$ space whereas  globally display a
$L^q$-decay, with $q\not=p$). The Wiener amalgam spaces,
introduced by H. Feichtinger in 1980 \cite{feichtinger80},
 enjoy this property and are the means  to perform such a finer analysis of local integrability and decay at infinity. For instance,  the desired Wiener amalgam space in
the  example above is denoted by $W(L^p,L^q)$; similarly,
one can consider other Banach spaces to measure the local behaviour of a function, e.g., $\Fur L^p$ instead of $L^p$, the related Wiener amalgam space being then $W(\Fur L^p,L^q)$,  and so on. In general, given
two (suitable) Banach function spaces $B$, $C$, the Wiener amalgam space
$W(B,C)$ is the space of functions which ``locally'' are in $B$
and ``globally'' in  $C$ (see \cite{feichtinger80} and Section 2
below for the precise definition).\par Here is a brief discussion
of the results proved in this paper. As usual we first establish
an estimate of dispersive type. Namely, in our setting, we prove
\begin{equation}\label{dispersive}
\|e^{it\Delta}u_0\|_{W(\Fur L^1, L^\infty)}\lesssim
\left(\frac{1+|t|}{t^{2}}\right)^{\frac{d}{2}}\|u_0\|_{W(\Fur
L^\infty,L^1)}.
\end{equation}
Notice that
$$L^1=W(L^1,L^1) \hookrightarrow W(\cF L^\infty,L^1)\,\,\,\, \mbox{and }\,\,\, W(\cF L^1,L^\infty)\hookrightarrow W(L^\infty,L^\infty)=L^\infty.$$
Comparing with the classical dispersive estimate \cite{Kato70}
\begin{equation}\label{classic}
\|e^{it\Delta}u_0\|_{L^\infty}\lesssim |t|^{-{d/2}}\|u_0\|_{L^1},
\end{equation}
we therefore get an improvement for every fixed $t\not=0$.
Indeed, we start from less regular
data $u_0\in W(\Fur L^\infty,L^1)$ (for example, a compactly supported Radon measure in
$\bR^d$), and we end up with a solution $u(t,\cdot)$ locally in $\Fur L^1$, which is strictly smaller than $L^\infty$.
Observe that in \eqref{dispersive} we recapture the classical time decay $|t|^{-d/2}$
as $|t|\to+\infty$, whereas  getting a better result locally in space costs  a worsening as $|t|\to0$: the factor $|t|^{-d/2}$ is replaced by $|t|^{-d}$ as $t\to 0$.

Next, we focus on space-time estimates.  Precisely, upon defining the Wiener
amalgam norms as
\begin{eqnarray}\label{norm} \|F\|_{W(L^{q_1},L^{q_2})_t W(\Fur
L^{r_1},L^{r_2})_x}:&=&\| \|F(t)\|_{W(\Fur
L^{r_1},L^{r_2})_x}\|_{W(L^{q_1},L^{q_2})_t}\\
&=&\|F\|_{W\left(L^{q_1}_tW(\Fur
L^{r_1}_x,L^{r_2}_x),L^{q_2}_t\right)}\nonumber
\end{eqnarray}
(the last equality shall be proved in Section 2), our result can
be stated as follows.

\begin{theorem}\label{prima}
Let $4<q,\tilde{q}\leq\infty$, $2\leq r,\tilde{r}\leq\infty$, such that
\begin{equation}\label{0000}\frac{2}{q}+\frac{d}{r}=\frac{d}{2},
\end{equation} and
similarly for $\tilde{q},\tilde{r}$. Then we have
the homogeneous Strichartz estimates
\begin{equation}\label{hom}\|e^{it\Delta} u_0\|_{W(L^{{q}/{2}},L^q)_t W(\Fur
L^{r^\prime},L^r)_x}\lesssim \|u_0\|_{L^2_x},
\end{equation}
the dual homogeneous Strichartz estimates
\begin{equation}\label{dh}
\|\int e^{-is\Delta} F(s)\,ds\|_{L^2}\lesssim
\|F\|_{W(L^{(\tilde{q}/{2})^\prime},L^{\tilde{q}^\prime})_t W(\Fur
L^{\tilde{r}},L^{\tilde{r}^\prime})_x},
\end{equation}
and the retarded Strichartz estimates
\begin{equation}\label{ret}
\|\int_{s<t} e^{i(t-s)\Delta} F(s)\,ds\|_{W(L^{q/2},L^q)_t W(\Fur
L^{r^\prime},L^r)_x}
\lesssim\|F\|_{W(L^{(\tilde{q}/{2})^\prime},L^{\tilde{q}^\prime})_t
W(\Fur L^{\tilde{r}},L^{\tilde{r}^\prime})_x}.
\end{equation}
\end{theorem}
\noindent The solution of \eqref{cp}, with $u_0\in L^2$, is therefore shown to be in the space \\ $L^{q/2}_{t,{\rm loc}}W(\Fur L^{r'},L^r)_x$ and, roughly speaking, its $W(\Fur L^{r'},L^r)_x$-norm has a $L^q_t$-decay at infinity, for $2/q+d/r=d/2$, $q>4$, $r\geq2$. Hence, for $q>4$, the result
is better than the classical one for large time, while locally we pay that
improvement: the classical $L^q_t$
 regularity is replaced by $L^{q/2}_t$.\par
In the case $(q,r)=(\infty,2)$ (included in Theorem \ref{prima}), we recapture the usual estimate
$$ \|e^{it\Delta} u_0\|_{L^\infty_t L^2_x}\lesssim \|u_0\|_{L^2_x}.$$

From a local point of view, the other cases
 are not comparable to those of \eqref{S1},
  since here we gain  in space (being $\Fur L^{r'}\subset L^r$ for $r\geq2$) but lose in time. In fact, it is important to observe that, for any given $s\geq1$, the inclusion $L^s \subset \left(\Fur L^{r^\prime}\right)_{\rm loc}$ is always {\it false} when $r>2$, so that estimate \eqref{hom} cannot be deduced from \eqref{S1} applied to the Schr\"odinger admissible pair $(q/2,s=2dq/(dq-8))$, $q>4$. In other terms, estimate \eqref{hom} contains information on the oscillations in $x$ of the solution, which cannot be extracted from \eqref{S1}.\par
 Observe that locally the $L^2_t$-norm would be obtained with $q=4$,
which corresponds to  our endpoint case. Namely,
let $P$ be the  endpoint
$$P:=(4,2d/(d-1)), \quad d>1;$$ then our version of the main result of \cite{keel} can be formulated as follows.
\begin{theorem}\label{endpoint}
For $(q,r)=P$, $d>1$, we have
\begin{equation}\label{hombis}
\|e^{it\Delta} u_0\|_{W(L^{2},L^4)_t W(\Fur
L^{r^\prime,2},L^r)_x}\lesssim \|u_0\|_{L^2_x},
\end{equation}
\begin{equation}\label{dhbis}
\|\int e^{-is\Delta} F(s)\,ds\|_{L^2}\lesssim
\|F\|_{W(L^{2},L^{4/3})_t W(\Fur
L^{r,2},L^{r^\prime})_x}.
\end{equation}
The retarded estimates \eqref{ret} still hold with $(q,r)$ satisfying \eqref{0000}, $q>4,r\geq2$, $(\tilde{q},\tilde{r})=P$, if one replaces $\Fur L^{\tilde{r}'}$ by $\Fur L^{\tilde{r}',2}$. Similarly it holds for $(q,r)=P$ and $(\tilde{q},\tilde{r})\not=P$ as above if one replaces $\Fur L^{r'}$ by $\Fur L^{r',2}$. It holds for both $(p,r)=(\tilde{p},\tilde{r})=P$ if one replaces $\Fur L^{r'}$ by $\Fur L^{r',2}$ and $\Fur L^{\tilde{r}'}$ by $\Fur L^{\tilde{r}',2}$.
\end{theorem}
Here $L^{r,2}$ is a Lorentz space (see \cite{steinweiss} and Section \ref{section2} below). We recall, $L^{r,2}\subset L^r$ and $L^{r'}\subset L^{r^\prime,2}$, for $r>2$. Indeed, here we attain slightly weaker estimates than in Theorem \ref{prima}.\\ Let us observe that, for $d=2$, there is no estimate of the form $\|e^{it\Delta}u_0\|_{L^2_t L^s_x}\lesssim\|u_0\|_{L^{2}}$, with $s\geq1$ (see \cite{tao2}). However, Theorem \ref{endpoint} above shows that the solution  $e^{it\Delta}u_0$ lies in the space $L^2_{t,{\rm loc}}(W(\Fur L^{4/3,2},L^4)_x)$ and that, roughly speaking, its $W(\Fur L^{4/3,2},L^4)_x$-norm has a $L^4_t$-decay at infinity.\par
A natural question  arises: can our results  be extended to the case $2\leq q<4$?
So, for instance, for $q=2$
one would obtain
$$\|e^{it\Delta} u_0\|_{W(L^{1},L^2)_t W(\Fur
L^{r^\prime},L^r)_x}\lesssim \|u_0\|_{L^2_x},\quad r=2d/(d-2).
$$
For sure this case cannot be treated with the techniques developed here, which use the Hardy-Littlewood-Sobolev's singular integral theory.

The Strichartz estimates of Theorems \ref{prima} and \ref{endpoint} can be applied, e.g., to the well-posedness of non-linear Schr\"odinger equations or of linear Schr\"odinger equations with time-dependent potentials. As an example, in Section \ref{section6} we combine our estimates with the methods of \cite{Dancona05} to deduce some estimates for the Schr\"odinger equation with a potential $V(t,x)\in L^\alpha_t L^p_x$.\par
We point out that interesting estimates for the operator $e^{it\Delta}$, for fixed $t$,  have been recently obtained in \cite{baoxiang} and \cite{benyi} using  the framework of modulation spaces. Such spaces are related to the Wiener amalgam spaces considered here via  Fourier transform. The  overlap with our results is Proposition \ref{chirp} below, which was first obtained there.\par
The paper is organized as follows. In Section \ref{section2} we
recall the definition and the main properties of the function
spaces used in this paper. In Section \ref{section3} we prove the
dispersive estimate and other fixed time estimates for the
solution of  \eqref{cp}. In Section \ref{section4} we prove Theorem
\ref{prima}. In Section \ref{section5} we
prove Theorem \ref{endpoint}. Finally in Section \ref{section6} we present the above mentioned application to Schr\"odinger equations with time-dependent potentials.\par
\par\medskip\noindent
\textbf{Notation.} We define $|x|^2=x\cdot x$, for $x\in\Ren$, where
$x\cdot y=xy$ is the scalar product on $\Ren$. The space of smooth
functions with compact support is denoted by $\cC_0^\infty(\rd)$,
the Schwartz class is  $\sch(\Ren)$, the space of tempered
distributions $\sch'(\Ren)$.    The Fourier transform is
normalized to be ${\hat
  {f}}(\o)=\Fur f(\o)=\int
f(t)e^{-2\pi i t\o}dt$.
 Translation and modulation operators ({\it time and frequency shifts}) are defined, respectively, by
$$ T_xf(t)=f(t-x)\quad{\rm and}\quad M_{\o}f(t)= e^{2\pi i \o
 t}f(t).$$
We have the formulas $(T_xf)\hat{} = M_{-x}{\hat {f}}$,  $(M_{\o}f)\hat{} =T_{\o}{\hat {f}}$, and $M_{\o}T_x=e^{2\pi i x\o}T_xM_{\o}$.
The notation $A\lesssim B$ means $A\leq c B$ for a suitable constant $c>0$, whereas $A
\asymp B$ means $c\inv A \leq B \leq c A$, for some $c\geq 1$. The
symbol $B_1 \hookrightarrow B_2$ denotes the continuous embedding of
the linear space $B_1$ into $B_2$.

\section{Function Spaces}\label{section2}
\noi \textbf{Lorentz spaces (\cite{stein93,steinweiss}).} We
recall that the Lorentz space $L^{p,q}$ on $\mathbb{R}^d$ is
defined as the space of measurable functions $f$ such that
\[
\|f\|^\ast_{pq}=\left(\frac{q}{p}\int_0^\infty [t^{1/p} f^\ast
(t)]^q\frac{dt}{t}\right)^{1/q}<\infty,
\]
when $1\leq p<\infty$, $1\leq q<\infty$, and
\[
\|f\|^\ast_{pq}=\sup_{t>0} t^{1/p}f^\ast(t)<\infty
\]
when $1\leq p\leq\infty$, $q=\infty$. Here, as usual,
$\lambda(s)=|\{|f|>s\}|$ denotes the distribution function of $f$
and $f^\ast(t)=\inf\{s:\lambda(s)\leq t\}$.\par We also recall
that the following equality holds:
$$\|f\|^\ast_{p\infty}=\sup_{s>0} s\lambda(s)^{1/p},
$$
which gives a characterization of $L^{p,\infty}$.\par 
One has
$L^{p,q_1}\subset L^{p,q_2}$ if $q_1\leq q_2$, and $L^{p,p}=L^p$.
Moreover, for $1<p<\infty$ and $1\leq q\leq\infty$, $L^{p,q}$ is a
normed space and its norm $\|\cdot\|_{L^{p,q}}$ is equivalent to the
above quasinorm $\|\cdot\|^\ast_{pq}$. \par The following
important result (\cite[Theorem 2]{triebel}, page 139) will be
crucial in the sequel. It generalizes the Hardy-Littlewood-Sobolev
fractional integration theorem (see e.g. \cite{stein}, page 119)
which corresponds to the model case of convolution by
$K(x)=|x|^{-\alpha}\in L^{d/\alpha,\infty}$, $0<\a<d$.

\begin{theorem}\label{convlor}
Let $1\leq p<q<\infty$, $0<\alpha<d$, with
$1/p=1/q+1-\alpha/d$. Then,
\begin{equation}\label{conv1}
L^p(\rd)\ast L^{d/\alpha,\infty}(\rd)\hookrightarrow L^q(\rd).
\end{equation}
\end{theorem}
 \par\medskip\noindent
\textbf{Wiener amalgam spaces  (\cite{feichtinger80,feichtinger83,
feichtinger90,fournier-stewart85,Fei98}).} Let $g \in
\cC_0^\infty(\bR^n)$ be a test function that satisfies
$\|g\|_{L^2}=1$. We will refer to $g$ as a window function.\\ Let
$B$ one of the following Banach spaces:
$L^p, \cF L^p$, $1\leq p\leq \infty$,  $L^{p,q}$, $1<p<\infty$, $1\leq q\leq \infty$, possibly valued in a Banach space, or also spaces obtained from these by real or complex interpolation.\\
Moreover, let $C$ be one of the following Banach spaces: $L^p$, $1\leq p\leq\infty$, or $L^{p,q}$,
 $1<p<\infty$, $1\leq q\leq \infty$, scalar valued.\\
For any given function $f$ which is locally in $B$ (i.e. $g
f\in B$, $\forall g\in\cC_0^\infty$), we fix $g\in\cC_0^\infty$ and  set $f_B(x)=\| fT_x
g\|_B$. Then, the {\it Wiener amalgam space} $W(B,C)$ with local
component $B$ and global component  $C$ is defined as the space of
all functions $f$ locally in $B$ such that $f_B\in C$. Endowed
with the norm
$\|f\|_{W(B,C)}:=\|f_B\|_C$, $W(B,C)$ is a Banach space. Besides, different choices of $g\in \cC_0^\infty$  generate the
same space and yield equivalent norms.
In particular, if we choose $B=\Fur L^1$ (the Fourier algebra), then the space of admissible windows for the Wiener amalgam spaces $W(\Fur L^1,C)$ can be enlarged to the
so-called Feichtinger algebra $W(\Fur L^1,L^1)$. Let us recall  that the Schwartz class $\sch$
  is dense in $W(\Fur L^1,L^1)$.\par
Observe that this definition mixes (\emph{amalgamates}) the local properties of functions in $B$ with the global properties of functions in $C$.

  Now, it is straightforward to prove the norm equality in \eqref{norm}.
  Precisely, for a fixed window function $g\in\cC_0^\infty(\R)$, we
  can write
  \begin{eqnarray*} \| \|F\|_{W(\Fur
L^{r_1},L^{r_2})_x}\|_{W(L^{q_1},L^{q_2})_t}&=&
\|\|T_ug(t)(\|F(t)\|_{W(\Fur
L^{r_1},L^{r_2})_x})\|_{L^{q_1}_t}\|_{L^{q_2}_u}\\
&=&\|\|\|\,|T_ug(t)|\,F(t)\|_{W(\Fur
L^{r_1},L^{r_2})_x}\|_{L^{q_1}_t}\|_{L^{q_2}_u}\\
&=& \|F\|_{W\left(L^{q_1}_t(W(\Fur
L^{r_1}_x,L^{r_2}_x)),L^{q_2}_t\right)}.
\end{eqnarray*}

  \par\medskip\noindent

  Hereafter we shall recall some useful  properties of the Wiener amalgam
   spaces.\\
  \begin{lemma}\label{WA}
  Let $B_i$, $C_i$, $i=1,2,3$ be Banach spaces  such that $W(B_i,C_i)$ are well
  defined. Then,
  \begin{itemize}
  \item[(i)] \emph{Convolution.}
  If $B_1\ast B_2\hookrightarrow B_3$ and $C_1\ast
  C_2\hookrightarrow C_3$, we have
  \begin{equation}\label{conv0}
  W(B_1,C_1)\ast W(B_2,C_2)\hookrightarrow W(B_3,C_3).
  \end{equation}
  \item[(ii)]\emph{Inclusions.} If $B_1 \hookrightarrow B_2$ and $C_1 \hookrightarrow C_2$,
   \begin{equation*}
   W(B_1,C_1)\hookrightarrow W(B_2,C_2).
  \end{equation*}
  \noindent Moreover, the inclusion of $B_1$ into $B_2$ need only hold ``locally'' and the inclusion of $C_1 $ into $C_2$  ``globally''.
   In particular, for $1\leq p_i,q_i\leq\infty$, $i=1,2$, we have
  \begin{equation}\label{lp}
  p_1\geq p_2\,\mbox{and}\,\, q_1\leq q_2\,\Longrightarrow W(L^{p_1},L^{q_1})\hookrightarrow
  W(L^{p_2},L^{q_2}).
  \end{equation}
  \item[(iii)]\emph{Complex interpolation.} For $0<\theta<1$, we
  have
\[
  [W(B_1,C_1),W(B_2,C_2)]_{[\theta]}=W\left([B_1,B_2]_{[\theta]},[C_1,C_2]_{[\theta]}\right),
  \]
if $C_1$ or $C_2$ has absolutely continuous norm.
    \item[(iv)] \emph{Duality.}
    If $B',C'$ are the topological dual spaces of the Banach spaces $B,C$ respectively, and
    the space of test functions $\cC_0^\infty$ is dense in both $B$ and $C$, then
\begin{equation}\label{duality}
W(B,C)'=W(B',C').
\end{equation}
  \end{itemize}
  \end{lemma}
  The proof of all these results can be found in
  (\cite{feichtinger80,feichtinger83,feichtinger90,Heil03}).\par
   Here  are instead some results on real interpolation theory we did not
    find explicitly established in the literature.
     We use the notation and terminology of \cite{triebel}.
\begin{proposition}
Let $\{A_0,A_1\}$ be an interpolation couple. For every $1\leq
p_0,p_1<\infty$, $0<\theta<1$, $1/p=(1-\theta)/p_0+\theta/p_1$ and
$p\leq q$ we have
\begin{equation}\label{inter1}
l^p\left((A_0,A_1)_{\theta,q}\right)\hookrightarrow
(l^{p_0}(A_0),l^{p_1}(A_1))_{\theta,q}.
\end{equation}
\end{proposition}
\begin{proof} Set $\eta=p\theta/p_1$ and $q=pr$, $r\geq1$. It follows from Theorem 1.4.2 of \cite{triebel}, page 29, that given $c=\{c_j\}\in l^p(A_0)+l^p(A_1)$ we have
\[
\|c\|_{(l^{p_0}(A_0),l^{p_1}(A_1))_{\theta,q}}^p\asymp\|t^{-\eta}\inf_{\substack{a+b=c\\
a\in l^{p_0}(A_0),\ b\in l^{p_1}(A_1)}}
\|a\|_{l^{p_0}(A_0)}^{p_0}+t\|b\|_{l^{p_1}(A_1)}^{p_1}\|_{L^r(\mathbb{R}_+,\frac{dt}{t})}.
\]
Hence,
\begin{align}
\|c\|_{(l^{p_0}(A_0),l^{p_1}(A_1))_{\theta,q}}^p&\asymp \|t^{-\eta} \inf_{\substack{a+b=c\\ a\in l^{p_0}(A_0),\ b\in l^{p_1}(A_1)}}\sum_{j\geq0} (\|a_j\|^{p_0}_{A_0}+t\|b_j\|^{p_1}_{A_1})\|_{{L^{r}(\mathbb{R}_+,\frac{dt}{t})}}\nonumber \\
&=\|t^{-\eta}\sum_{j\geq0}\inf_{\substack{a_j+b_j=c_j\\ a_j\in A_0
,\ b_j\in A_1}}
(\|a_j\|^{p_0}_{A_0}+t\|b_j\|^{p_1}_{A_1})\|_{{L^{r}(\mathbb{R}_+,\frac{dt}{t})}}.
\end{align}
By Minkowski's inequality we deduce,
\begin{align}
\|c\|_{(l^{p_0}(A_0),l^{p_1}(A_1))_{\theta,q}}^p&\lesssim\sum_{j\geq0} \|t^{-\eta}\inf_{\substack{a_j+b_j=c_j\\ a_j\in A_0 ,\ b_j\in A_1}} (\|a_j\|^{p_0}_{A_0}+t\|b_j\|^{p_1}_{A_1})\|_{{L^{r}(\mathbb{R}_+,\frac{dt}{t})}}\nonumber\\
&\asymp\|c\|^p_{l^p\left((A_0,A_1)_{\theta,q}\right)}.
\end{align}
\end{proof}

Consider now a partition of unity \footnote{Such a partition of
unity can be constructed as follows. Take $\chi\in
\mathcal{C}^\infty_0(\mathbb{R}^d)$, $0\leq\chi\leq1$, $\chi=1$ on
$[-1,1]^d$, $\chi=0$ away from $[-2,2]^d$. Set
$\Phi(x)=\sum_{\alpha\in\mathbb{Z}^d} \chi( x-\alpha)$. Since the
sum is locally finite, $\Phi$ is well defined and smooth. Moreover
$\Phi(x+\beta)=\Phi(x)$ for every $\beta\in\mathbb{Z}^d$, and also
$\Phi\geq 1$.  Hence it suffices to take
$\phi_\alpha=T_\alpha\chi/\Phi=T_\alpha(\chi/\Phi)$} given by
functions $\phi_\alpha\in \mathcal{C}^\infty_0(\mathbb{R}^d)$,
$\alpha\in\mathbb{Z}^d$, with $\phi_\alpha=T_\alpha \phi$, ${\rm
supp\, \phi}\subset[-2,2]^d$. Thus ${\rm supp}\,\phi_\alpha\subset
\alpha+[-2,2]^d$. Let then $\psi\in \mathcal{C}^\infty_0(\mathbb{R}^d)$,
$\psi=1$ on $[-2,2]^d$ and $\psi=0$ away from $[-4,4]^d$, and set
$\psi_\alpha=T_\alpha\psi$. Observe that there is a constant $C_d$
such that,
\begin{equation}\label{m1}
\forall\alpha\in\mathbb{Z}^d,\  \#\{\beta\in\mathbb{Z}^d: {\rm
supp}\, \psi_\beta\cap {\rm supp}\,\phi_\alpha\not=\emptyset\}\leq
C_d
\end{equation}
and
\begin{equation}\label{m2}
\forall\alpha\in\mathbb{Z}^d,\  \#\{\beta\in\mathbb{Z}^d: {\rm
supp}\, \psi_\alpha\cap {\rm supp}\,\phi_\beta\not=\emptyset\}\leq
C_d.
\end{equation}
The linear operators
\[
S: L^1_{\rm loc}\to(L^1_{\rm loc})^{\mathbb{Z}^d},\quad R:
(L^1_{\rm loc})^{\mathbb{Z}^d}\to L^1_{\rm loc},
\]
defined  by
\begin{equation}\label{retr}
Sf=\{f\phi_\alpha\}_\alpha,\quad
R(\{u_\alpha\}_\alpha)=\sum_{\alpha} u_\alpha{\psi}_\alpha,
\end{equation}
enjoy the following properties (see \cite[Remark 2.2]{feichtinger83}).
\begin{proposition}\label{inter2} 
We have $RS={\rm Id}$ on $L^1_{\rm loc}$ and, for every local
component $B$, as at the beginning of this section, and every
$p\geq1$, we have
\begin{equation}\label{n1}
S: W(B,L^p)\to l^p(B),
\end{equation}
and
\begin{equation}\label{n2}
R: l^p(B)\to W(B,L^p)
\end{equation}
continuously.
\end{proposition}
\begin{proof} The equality $RS={\rm Id}$ on $L^1_{\rm loc}$ is clear, whereas
\eqref{n1} follows at once from Remark 4 of \cite{feichtinger80}.
To prove \eqref{n2} we observe that
\begin{align}
\|R(\{u_\alpha\}_\alpha)\|^p_{W(B,L^p)}&\asymp\sum_\beta\|\sum_\alpha u_\alpha \psi_\alpha\phi_\beta\|^p_B,\quad (\textrm{by Remark 4 of \cite{feichtinger80}}) \\
&\lesssim\sum_\beta\sum_\alpha\|u_\alpha\psi_\alpha\phi_\beta\|^p_B\quad (\textrm{by \eqref{m1}})\\
&\lesssim \sum_\alpha\sum_\beta\|\psi_\alpha\phi_\beta\|^p_{\Fur L^1} \|u_\alpha\|^p_B\\
&\lesssim \|\{u_\alpha\}_\alpha\|^p_{l^p(B)}.
\end{align}
In the last inequality we used the fact that, by \eqref{m2},
\[
\sum_\beta\|\psi_\alpha\phi_\beta\|^p_{\Fur L^1}\lesssim
\sum_{\beta:\ {\rm supp}\, \psi_\alpha\cap {\rm
supp}\,\phi_\beta\not=\emptyset} \|\psi_\alpha\|^p_{\Fur L^1}
\|\phi_\beta\|^p_{\Fur L^1}\leq C_d \|\psi_0\|^p_{\Fur L^1}
\|\phi_0\|^p_{\Fur L^1}.
\]
\end{proof}
\begin{proposition}\label{inter9}
Given two local components $B_0,B_1$ as at the beginning of this
section, for every $1\leq p_0,p_1<\infty$, $0<\theta<1$,
$1/p=(1-\theta)/p_0+\theta/p_1$, and $p\leq q$ we have
\[
W\left((B_0,B_1)_{\theta,q},L^p\right)\hookrightarrow\left(W(B_0,L^{p_0}),W(B_1,L^{p_1})\right)_{\theta,q}.
\]
\end{proposition}
\begin{proof}
Let $R$ and $S$ be the operators defined above. Then, by
Proposition \ref{inter2},
\begin{align}
\|f\| _{\left(W(B_0,L^{p_0}),W(B_1,L^{p_1})\right)_{\theta,q}}&=\|RS f\|_{\left(W(B_0,L^{p_0}),W(B_1,L^{p_1})\right)_{\theta,q}}\nonumber\\
&\lesssim\|Sf\|_{(l^{p_0}(B_0),l^{p_1}(B_1))_{\theta,q}}\nonumber\\
&\lesssim \|Sf\|_{l^p\left((B_0,B_1)_{\theta,q}\right)},
\label{terz}\\
&\lesssim \|f\|_{W\left((B_0,B_1)_{\theta,q},L^p\right)},\nonumber
\end{align}
where for \eqref{terz} we used \eqref{inter1}.
This concludes the proof.
\end{proof}

\section{Fixed time estimates}\label{section3}
In this section we study estimates for the solution $u(t,x)$ of the Cauchy problem \eqref{cp}, for fixed $t$. We take advantage of the explicit formula for the solution
\begin{equation}\label{sol}
u(t,x)=(K_t\ast u_0)(x)
\end{equation}
where
\begin{equation}\label{chirp0}
K_t(x)=\frac{1}{(4\pi i t)^{d/2}}e^{i|x|^2/(4t)}.
\end{equation}
Precisely, we are going to show that the function in
\eqref{chirp0} is in the Wiener amalgam space $W(\cF L^1,
L^\infty)$ and we shall compute its norm. This goal is attained
thanks to the nice choice of the Gaussian function $e^{-\pi
|x|^2}\in \cS(\rd)$, as a fitting window function. Then we will make
use of the convolution properties of Wiener amalgam spaces.\par We
first remind the \ft \,of the Gaussian function (see, e.g.,
\cite{folland89}).
\begin{lemma}
For ${\rm Re}\,c\geq 0$, $c\not=0$. Let $\f_c(x)=e^{-\pi |x|^2/c}$, then
\begin{equation}\label{ftGauss}
\widehat{\f_c}(\o)=c^{d/2}\f_{1/c}(\o),\quad \o\in\rd,
\end{equation}
where the square root is chosen to have a positive real part.
\end{lemma}

\begin{proposition}\label{chirp} (cf. \cite{baoxiang,benyi})
For $a\in\mathbb{R}$, $a\not=0$, let $f_a(x)=(a i)^{-d/2}e^{-\pi |x|^2/(ai)}$. Then
$f_a\in W(\cF L^1, L^\infty)$, with
\begin{equation}\label{chirpnorm}
\|f_a\|_{W(\cF L^1,
L^\infty)}=\left(\frac{1+a^2}{a^{4}}\right)^{\frac{d}{4}}.\end{equation}
\end{proposition}
\begin{proof}
By definition (see, e.g., \cite{Fei98,Heil03}),
\begin{equation*}
\|f_a\|_{W(\cF L^1, L^\infty)}=\sup_{x\in\rd}\|f_aT_x {g}\|_{\cF
L^1},
\end{equation*}
for some  non-zero window $g\in W(\cF
L^1,L^1)$ (different windows give equivalent norms). We then
choose $g=e^{-\pi |x|^2}$, and, observing that $\hat{g}=g$, we can
write
$$\|f_aT_x{g}\|_{\cF L^1}=\|\widehat{f_aT_x {g}}\|_{L^1}=\|\hat{f_a}\ast M_{-x}g \|_{L^1}.
$$
Using \eqref{ftGauss} with   $c=ai $ we compute the \ft\, of $f$
that reveals to be
$$\hat{f_a}(\o)=(ai )^{-d/2}(ai)^{d/2}e^{-\pi ai\o^2} =e^{-\pi a i\o^2}.$$
Thereby,
\begin{align*}
\|\hat{f_a}\ast M_{-x}g \|_{L^1}&=\intrd\left|\intrd e^{-\pi a i(\o-y)^2}e^{-2\pi i x y}e^{-\pi^2 |y|^2}\,dy\right|d\o \\
&=\intrd\left|e^{-\pi a i\o^2}\intrd e^{-2\pi i(x-a\o)y}e^{-\pi(1+a i )|y|^2}\,dy\right|d\o\\
&=\intrd\left|\cF (e^{-\pi(1+a i )|y|^2})(x-a\o)\right|d\o\\
&=\intrd\left|(1+a i)^{-d/2} e^{-\pi(x-a\o)^2/(1+a i )}\right|d\o,
\end{align*}
 where in the last equality we use \eqref{ftGauss}. Performing the change of variables $ x-a\o=z$, hence $d\o=a^{-d}dz $ and
 observing that $|(1+a i)^{-d/2}|=(1+a^2)^{-d/4}$, we can write
  \begin{align*}
  \|\hat{f_a}\ast M_{-x}g \|_{L^1}&=(1+a^2)^{-d/4}a ^{-d}\intrd \left| e^{-\pi |z|^2/(1+a i )} dz\right|\\
  &=(1+a^2)^{-d/4}a ^{-d}\intrd e^{-\pi |z|^2/(1+a^2)} dz\\
  &=(1+a^2)^{-d/4}a ^{-d}(1+a^2)^{d/2}\\
  &=\left(\frac{1+a^2 }{a^4 }\right)^{d/4}.
  \end{align*}
  Since the right-hand side does not depend on $x$, taking the supremum on $\rd$ with respect to the $x$-variable we attain the desired estimate.
\end{proof}
\begin{lemma}
It turns out
\begin{equation}\label{eq01}
W(\Fur L^1,L^\infty)\ast W(\Fur L^\infty,L^1)\hookrightarrow   W(\Fur L^1,L^\infty).
\end{equation}
\end{lemma}
\begin{proof}
This is a consequence of the convolution relations for Wiener
amalgam spaces in Lemma \ref{WA} (i), being $\Fur L^1\ast \Fur
L^\infty= \Fur (L^1\cdot L^\infty)=\Fur L^1$ and $L^\infty\ast L^1
\hookrightarrow L^\infty$.
\end{proof}
\begin{prop}
We have
\begin{equation}\label{est1}
\|e^{it\Delta}u_0\|_{W(\Fur L^1, L^\infty)}\lesssim
\left(\frac{1+|t|}{t^{2}}\right)^{\frac{d}{2}}\|u_0\|_{W(\Fur
L^\infty,L^1)}.
\end{equation}
\end{prop}
\begin{proof} We use the explicit representation of the Schr\"odinger evolution operator $ e^{it\Delta}u_0(x)=(K_t\ast
u_0)(x)$. From \eqref{chirpnorm} (with $a=4\pi t$) we infer
$$ \|K_t\|_{W(\Fur L^1, L^\infty)} \asymp
\left(\frac{1+|t|}{t^{2}}\right)^{\frac{d}{2}}.$$ Finally, the
convolution relations  \eqref{eq01} yield the desired result.
\end{proof}
\begin{theorem}
For $2\leq r\leq \infty$ we have
\begin{equation}\label{est2}
\|e^{it\Delta}u_0\|_{W(\Fur L^{r^\prime}, L^r)}\lesssim \left(\frac{1+|t|}{t^2}\right)^{d\left(\frac{1}{2}-\frac{1}{r}\right)}\|u_0\|_{W(\Fur L^r,L^{r^\prime})}.
\end{equation}
\end{theorem}
\begin{proof}
Estimate \eqref{est2} follows by interpolating \eqref{est1} with the $L^2$ conservation law
\begin{equation}\label{l2}
\|e^{it\Delta}u_0\|_{L^2}=\|u_0\|_{L^2}.
\end{equation}
Indeed, $L^2=W(\Fur L^2,L^2)$. By  Lemma \ref{WA}, item (iii),  for
$0<\theta=2/r<1$,
\[
\left[W(\Fur L^1, L^\infty),W(\Fur L^2,L^2)\right]_{[\theta]}=W\left([\Fur L^1,\Fur L^2]_{[\theta]}, [L^\infty,L^2]_{[\theta]}\right)=W(\Fur L^{r^\prime}, L^r)
\]
and
\[
\left[W(\Fur L^\infty, L^1),W(\Fur L^2,L^2)\right]_{[\theta]}=W\left([\Fur L^\infty,\Fur L^2]_{[\theta]}, [L^1,L^2]_{[\theta]}\right)=W(\Fur L^{r}, L^{r^\prime}).
\]
\end{proof}
\begin{remark}\rm
As well known (see e.g. (2.23) of \cite{tao}), the $L^p$
fixed time estimates for the solution of \eqref{cp} read
\begin{equation}\label{est3}
\|e^{it\Delta}u_0\|_{L^{r}(\bR^d)}\lesssim |t|^{-d\left(\frac{1}{2}-\frac{1}{r}\right)}\|u_0\|_{L^{r^\prime}(\bR^d)},\quad 2\leq r\leq\infty.
\end{equation}
For $2\leq r\leq\infty$, $\Fur L^{r^\prime}\hookrightarrow L^r$, and the inclusion relations for Wiener amalgam spaces in Lemma \ref{WA}, item (ii), yield
 $W(\Fur L^{r^\prime},L^r)\hookrightarrow W(L^r,L^r)= L^r$ and $L^{r^\prime}=W(L^{r^\prime},L^{r^\prime})\hookrightarrow W(\Fur L^{r},L^{r^\prime})$. Thereby the estimates \eqref{est2} are indeed an improvement of \eqref{est3} for every fixed $t\not=0$, and also uniformly for $|t|>c>0$.
\end{remark}

\section{Strichartz estimates: non-endpoint case}\label{section4}
In this section we shall prove Theorem \ref{prima}. To this aim, we need some preliminary results.  For $0<\alpha<1/2$, let $\phi_\a$ be
the non-negative real function defined by
\begin{equation}\label{phialpha}
\phi_\alpha(t)=|t|^{-\alpha}+|t|^{-2\alpha},\quad t\in\R,\,t\not=0.
\end{equation}
Of course, $\phi_\a\in L^1_{loc}(\bR)$. The next lemma sets $\phi_\a$ in a suitable Wiener amalgam space.
\begin{lemma}\label{lemma1} We have
\[
\phi_\alpha\in W(L^{1/(2\alpha),\infty},L^{1/\alpha,\infty}).
\]
\end{lemma}
\begin{proof}
Let
\[
G(x)=\|\phi_\alpha \chi_{[x-1,x+1]}\|_{L^{1/(2\alpha),\infty}},
\]
then, it suffices to prove that
\begin{equation}\label{claim}
G(x)\leq \frac{C}{|x|^\alpha}.
\end{equation}
We can suppose $x>0$ ($G(x)$ is an even function). Now, if
$0<x\leq 2$, we have
\[
G(x)\leq \|\phi_\alpha\chi_{[-1,3]}\|_{L^{1/(2\alpha),\infty}}.
\]
On the other hand, when $x>2$, we use the fact that $\phi_\alpha(t)\leq 2|t|^{-\alpha}$ for $|t|\geq1$. Hence the distribution function $\lambda(s)$ of $\phi_\alpha(t)\chi_{[x-1,x+1]}(t)$ satisfies the estimate
\[
\lambda(s)\leq
\begin{cases}
2&{\rm if}\ s<2(x+1)^{-\alpha}\\
2^{1/\alpha}s^{-1/\alpha}-(x-1)&{\rm if}\ 2(x+1)^{-\alpha}\leq s\leq 2(x-1)^{-\alpha}\\
0& {\rm if}\ s>2(x-1)^{-\alpha}.
\end{cases}
\]
As a consequence,
\[
G(x)\asymp \sup_{s>0}\{s\lambda(s)^{2\alpha}\}\leq 2^{2\alpha+1}(x+1)^{-\alpha}.
\]
This proves \eqref{claim}.
\end{proof}
\begin{lemma}\label{0est2}
For $0<\alpha<1/2$, let $\phi_\alpha$ be the function defined in
\eqref{phialpha}.  Then,
\begin{equation}\label{conv2}
\|F\ast \phi_\alpha\|_{W(L^{1/\alpha},L^{2/\alpha})}\lesssim
\|F\|_{W(L^{{(1/\alpha)}^\prime},L^{(2/\alpha)^\prime})}.
\end{equation}
\end{lemma}
\begin{proof}
Lemma \ref{lemma1} above shows that, locally, $\phi_\alpha\in
L^{1/(2\a),\infty}$. Using the fractional integration Theorem
\ref{convlor} we infer
$$ L^p\ast L^{1/(2\a),\infty}\,\hookrightarrow\, L^s,\quad \mbox{with}\quad\frac 1p=\frac1s+1-2\a;$$
if we set $p=s'$, then $s=1/\a$ and $L^{(1/\alpha)^\prime}\ast
L^{1/(2\a),\infty} \hookrightarrow L^{1/\alpha}$. Globally,  $\phi_\alpha\in L^{1/\a,\infty}$ and the same argument
 gives
 $ L^{(2/\alpha)^\prime}\ast L^{1/\a,\infty}\hookrightarrow L^{2/\a}.$

Finally, the convolution relations for Wiener amalgam spaces
\eqref{conv0} glue together the local and global properties and
provide \eqref{conv2}.
\end{proof}
\begin{remark}\rm
Condition $\a<1/2$ is necessary if we want at least $\phi_\a\in
L^1_{loc}$ and solutions with local time estimate  in $L^p$ spaces
rather than rougher spaces, subsets of $\cS'$. This constraint will yield the threshold
$q=4$ in the Strichartz estimates.
\end{remark}

\begin{proof}[Proof of Theorem \ref{prima}]
We first  prove the estimate \eqref{hom}.

The case: $q=\infty$, $r=2$ follows at once from the conservation
law \eqref{l2}. Indeed, $W(L^{\infty},L^\infty)_t=L^\infty_t$ and
$ W(\Fur L^{2},L^2)_x=L^2_x$, so that,  taking the supremum over
$t$ in $\|e^{it\Delta} u_0\|_{L^2_x}=\|u_0\|_{L^2_x}$,  we attain
the claim.

 To prove the remaining cases, we can apply the usual
$TT^\ast$ method (or ``orthogonality principle", see \cite[Lemma 2.1]{GinibreVelo92} or \cite[page 353]{stein}), because of the H\"older's type inequality
\begin{equation}\label{holder}
|\langle F,G\rangle_{L^2_t L^2_x}|\leq \|F\|_{W(L^s,L^q)_t W(\Fur
L^{r^\prime},L^r))_x}\|G\|_{W(L^{s^\prime},L^{q^\prime})_t W(\Fur
L^{r},L^{r^\prime})_x},
\end{equation}
which can be proved directly from the definition of these spaces.\\
As a consequence, it suffices to prove the estimate
\begin{equation}\label{aus}
\|\int e^{i(t-s)\Delta}F(s)\,ds\|_{W({L^{q/2}},{L^q})_t W(\Fur
L^{r^\prime},L^r)_x}\lesssim\|F\|_{W({L^{\left(q/2\right)^\prime}},L^{q^\prime})_t
W(\Fur L^r,L^{r^\prime})_x}.
\end{equation}
Now, set $\alpha=d(1/2-1/r)=2/q$. Then, by
\eqref{est2} and Lemma \ref{0est2},
\begin{align*}
\|\int e^{i(t-s)\Delta}F(s)\,ds&\|_{W(L^{q/2},{L^q})_t W(\Fur
L^{r^\prime},L^r)_x}\\
&\leq \|\int\|e^{i(t-s)\Delta} F(s)\|_{W(\Fur
L^{r^\prime},L^r)_x}\,ds\|_{W({L^{q/2}},{L^q})_t}\\
&\lesssim\|\|F(t)\|_{W(\Fur L^{r^\prime},L^r)_x}\ast
\phi_\alpha(t)\|_{W(L^{q/2},{L^q})_t}\\
&\lesssim \|F\|_{W(L^{\left(q/2\right)^\prime},L^{q^\prime})_t
W(\Fur L^r,L^{r^\prime})_x}.
\end{align*}
The estimate \eqref{dh} follows  from \eqref{hom} by duality.
\par
Consider now the retarded Strichartz estimate \eqref{ret}. By complex interpolation (Lemma \ref{WA}, (iii)), in order to get \eqref{ret} with $(1/q,1/r)$, $(1/\tilde{q},1/\tilde{r})$, $(1/\infty,1/2)$ collinear it suffices to prove \eqref{ret} in the three
cases $(\tilde{q},\tilde{r})=(q,r)$, $(q,r)=(\infty,2)$,
$(\tilde{q},\tilde{r})=(\infty,2)$, as shown in Figure 1.\\
\vspace{-0.2cm}
 \begin{center}
           \includegraphics{figstr1.1}
            \\
           $ $
\end{center}
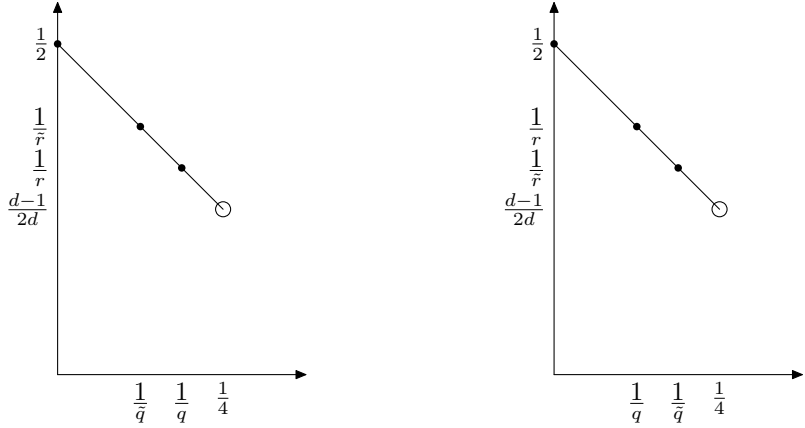
           \begin{center}{ Figure 1. The complex interpolation method is applied in these two  cases. }
           \end{center}
  \vspace{1.2cm}
Now, the case $(\tilde{q},\tilde{r})=(q,r)$ follows from
\eqref{aus} with $\chi_{\{s<t\}}F$ in place of $F$.\\
The case $(q,r)=(\infty,2)$ follows because
\begin{align*}
\|\int_{s<t} e^{i(t-s)\Delta} F(s)\,ds\|_{L^2_x}
&=\|\int_{s<t} e^{-is\Delta} F(s)\,ds\|_{L^2_x}\\
&\lesssim
\|F\|_{W(L^{(\tilde{q}/2)^\prime},L^{\tilde{q}^\prime})_t
W(\Fur L^{\tilde{r}},L^{\tilde{r}^\prime})_x},
\end{align*}
where we applied \eqref{l2} and then $\eqref{dh}$ with
$\chi_{\{s<t\}}F$ in place of $F$.\\
Finally, this latter argument applied to the adjoint operator
\[
G\mapsto \int _{t>s} e^{-i(t-s)\Delta} G(t)\,dt
\]
gives the case $(\tilde{q},\tilde{r})=(\infty,2)$.
\end{proof}
\section{Strichartz estimates: endpoint case}\label{section5}
In this section we prove the estimates in Theorem \ref{endpoint}. Hence,
$q=4,\,r=\frac{2d}{d-1}$, or $\tilde{q}=4,\,\tilde{r}=\frac{2d}{d-1}$, $d>1$.
We follow the pattern in Keel-Tao \cite{keel}. Indeed,
we study bilinear form estimates rather than operator
estimates. This is achieved via a dyadic decomposition in time (see \eqref{dec}), and estimates of each dyadic contribution (see \eqref{primolemma} and \eqref{lammasecondo}). Finally we conclude by a lemma of real interpolation theory. The proof will require however some different technical issues, due to the different nature of the Wiener amalgam spaces.\par
First we prove \eqref{hombis} and \eqref{dhbis}. Let therefore $$r=\frac{2d}{d-1},\quad d>1.$$ By the same duality arguments as the ones used in the
previous section, we observe that it suffices to prove
\eqref{hom}. This is equivalent to the bilinear estimate
\[
|\iint \langle e^{-is\Delta} F(s),
e^{-it\Delta}G(t)\rangle\,ds\,dt|\lesssim \|F\|_{W(L^2,L^{4/3})_t
W(\Fur L^{{r,2}},L^{{r}^\prime})_x} \|G\|_{W(L^2,L^{4/3})_t W(\Fur
L^{{r,2}},L^{{r}^\prime})_x}.
\]
By symmetry, it is enough to prove
\begin{equation}\label{aim}
|T(F,G)|\lesssim \|F\|_{W(L^2,L^{4/3})_t W(\Fur
L^{{r},2},L^{{r}^\prime})_x} \|G\|_{W(L^2,L^{4/3})_t W(\Fur
L^{{r},2},L^{{r}^\prime})_x},
\end{equation}
where
\[
T(F,G)=\iint_{s<t} \langle e^{-is\Delta} F(s),
e^{-it\Delta}G(t)\rangle\,ds\,dt.
\]
The form  $T(F,G)$ can be decomposed dyadically
as
\begin{equation}\label{dec} T=\tilde{T}+\sum_{j\leq 0} T_j,
\end{equation}
with
\[
\tilde{T}(F,G)=\iint_{s\leq t-2} \langle e^{-is\Delta} F(s),
e^{-it\Delta}G(t)\rangle\,ds\,dt
\]
and
\begin{equation}\label{tj}
T_j(F,G)=\iint_{t-2^{j+1}<s\leq t-2^j} \langle e^{-is\Delta} F(s),
e^{-it\Delta}G(t)\rangle\,ds\,dt.
\end{equation}
In the sequel we shall study the behaviour of $\tilde{T}$ and $T_j$
separately.
\begin{lemma}
We have
\begin{equation}\label{primolemma}
|\tilde{T}(F,G)|\lesssim \|F\|_{W(L^2,L^{4/3})_t W(\Fur
L^{{r}},L^{{r}^\prime})_x} \|G\|_{W(L^2,L^{4/3})_t W(\Fur
L^{{r}},L^{{r}^\prime})_x},
\end{equation}
\end{lemma}
\begin{proof}
It follows from the duality properties \eqref{duality} and the
space estimate \eqref{est2} that
\begin{eqnarray*}
|\langle e^{-is\Delta} F(s), e^{-it\Delta}G(t)\rangle|&=&
|\langle  F(s), e^{i(s-t)\Delta}G(t)\rangle|\\
&\lesssim&\|F(s)\|_{W(\Fur L^{{r}},L^{{r}^\prime})_x}
\|e^{i(s-t)\Delta}G(t) \|_{W(\Fur
L^{{r}^\prime},L^{{r}})_x}\\
&\lesssim& \|F(s)\|_{W(\Fur L^{{r}},L^{{r}^\prime})_x}\!\!\left(
\frac{1+|s-t|}{(s-t)^2}\right)^{d\left(\frac12-\frac1r\right)}\!\!\!
\|G(t)\|_{W(\Fur L^{{r}},L^{{r}^\prime})_x}.
\end{eqnarray*}
Since $d/(1/2-1/r)=1/2$, for $|s-t|\geq 2$ we have
$$\left(
\frac{1+|s-t|}{(s-t)^2}\right)^{d\left(\frac12-\frac1r\right)}\lesssim(1+|s-t|)^{-\frac12}.
$$
Thus, the form $\tilde{T}(F,G)$ can be controlled by the following
majorizations:
\begin{eqnarray*}
|\tilde{T}(F,G)|\!\!\!\!\!&&\lesssim\iint_{s\leq
t-2}\|F(s)\|_{W(\Fur
L^{{r}},L^{{r}^\prime})_x}(1+|s-t|)^{-\frac12}\|G(t)\|_{W(\Fur
L^{{r}},L^{{r}^\prime})_x}\,dsdt\\
&&\lesssim\|(1+|t|)^{-\frac12}\ast\|F(t)\|_{W(\Fur
L^{{r}},L^{{r}^\prime})_x}\|_{W(L^2,L^{(4/3)^\prime})_t}\\
&&\qquad\qquad\qquad\qquad\qquad\qquad\qquad\cdot\,\,\|\,\|G(t)\|_{W(\Fur
L^{{r}},L^{{r}^\prime})_x}\|_{W(L^2,L^{4/3})_t};
\end{eqnarray*}
in the last estimate we used the duality property for
 Wiener amalgam spaces \eqref{duality}. Notice that
 $(4/3)^\prime=4$.

  Next, the convolution relation  $\eqref{conv1}$, for $p=4/3$,
 $d=1$ and $\a=1/2$, reads $L^{4/3}\ast
 L^{2,\infty}\hookrightarrow L^4$, whereas Young's inequality
 gives $L^1\ast L^2\hookrightarrow L^2$. Applying the former
 estimate globally and the latter locally, we use the convolution relations for Wiener amalgam spaces \eqref{conv0} and infer
 $$ W(L^1,L^{2,\infty})\ast W(L^2,L^{4/3})\hookrightarrow W(L^2, L^4).$$
Since $4/3<4$ the inclusion relations \eqref{lp} give
$W(L^2,L^{4/3})\hookrightarrow W(L^2, L^4)$. This argument ends
our proof; indeed, it is straightforward to see that
$(1+|t|)^{-\frac12}\in
 W(L^1,L^{2,\infty})$, so that
\begin{eqnarray*}|\tilde{T}(F,G)|\!\!\!\!\!&&\lesssim\|(1+|t|)^{-\frac12}\ast\|F(t)\|_{W(\Fur
L^{{r}},L^{{r}^\prime})_x}\|_{W(L^2,L^4)_t}\|\,\|G(t)\|_{W(\Fur
L^{{r}},L^{{r}^\prime})_x}\|_{W(L^2,L^{4/3})_t}\\
&&\lesssim \|F\|_{W(L^2,L^{4/3})_tW(\Fur
L^{{r}},L^{{r}^\prime})_x}\|G\|_{W(L^2,L^{4/3})_tW(\Fur
L^{{r}},L^{{r}^\prime})_x},
\end{eqnarray*}
 as desired.
\end{proof}\\
For $a,b\geq 1$, we define
\[
\beta(a,b)=d-1-\frac{d}{a}-\frac{d}{b}.
\]
\begin{lemma}
\begin{equation}\label{lammasecondo}
|T_j(F,G)|\lesssim 2^{-j\beta(a,b)}\|F\|_{W(L^2,L^{4/3})_t W(\Fur L^a,L^{a^\prime})}\|G\|_{W(L^2,L^{4/3})_t W(\Fur L^b,L^{b^\prime})},
\end{equation}
for $(1/a,1/b)$ in a neighborhood of $(1/r,1/r)$.
\end{lemma}
\begin{proof}
Observe that here $r<\infty$. Then, the  result follows by
complex interpolation (Lemma \ref{WA}, (iii)) from the following
cases:
\begin{itemize}
\item[(i)] $a=\infty$, $b=\infty$, \item[(ii)] $2\leq a<r$, $b=2$,
\item[(iii)] $a=2$, $2\leq b<r$.
\end{itemize}
{\it Case (i)}. We need to show the estimate
\begin{equation}\label{(i)}|T_j(F,G)|\lesssim 2^{-j(d-1)} \|F\|_{W(L^2,L^{4/3})_t W(\Fur
L^\infty,L^1)}\|G\|_{W(L^2,L^{4/3})_t W(\Fur L^\infty,L^1)}.
\end{equation}
By \eqref{est1} and \eqref{phialpha} we have
\[
|\langle e^{-is\Delta} F(s), e^{-it\Delta} G(t)\rangle|\lesssim \phi_{d/2}(t-s)\|F(s)\|_{W(\Fur L^\infty,L^1)}\|G(t)\|_{W(\Fur L^\infty,L^1)}.
\]
Thereby
\[
|T_j(F,G)|\lesssim \phi_{d/2}(2^j)\|F\|_{L^1_t(W(\Fur L^\infty,L^1)_x)}\|G\|_{L^1_t(W(\Fur L^\infty,L^1)_x)}.
\]
We can of course assume $F$ and $G$ compactly supported, with respect to the time, in intervals of duration $\asymp 2^j$. As a consequence of H\"older's inequality,
\[
\|F(s)\|_{L^1_t(W(\Fur L^\infty,L^1)_x)}\lesssim 2^{j/2} \|F\|_{L^2_t(W(\Fur L^\infty,L^1)_x)},
\]
and similarly for $G$. Hence
\[
|T_j(F,G)|\lesssim \phi_{d/2}(2^j)2^j\|F\|_{L^2_t(W(\Fur L^\infty,L^1)_x)}\|G\|_{L^2_t(W(\Fur L^\infty,L^1)_x)}.
\]
Since $\phi_{d/2}(2^j)\lesssim 2^{-dj}$ (recall, $j\leq0$) and
\[
W(L^2_t(W(\Fur L^\infty,L^1)_x), L^{4/3})\hookrightarrow W(L^2_t(W(\Fur L^\infty,L^1)_x),L^{2}_t)=L^2_t(W(\Fur L^\infty,L^1)_x),
\]
we attain  the desired  estimate \eqref{(i)}.\\
 {\it Case (ii)}. We have to show
\begin{equation}\label{(ii)}|T_j(F,G)|\lesssim 2^{-j(\frac d2-1-\frac d a)} \|F\|_{W(L^2,L^{4/3})_t W(\Fur
L^a,L^{a^\prime})_x}\|G\|_{W(L^2,L^{4/3})_t  L^2_x}.
\end{equation}
Using similar arguments to the previous case we obtain
\begin{equation}\label{unozero}
|T_j(F,G)|\lesssim \sup_{t}\|\int_{t-2^{j+1}<s\leq t-2^j}
e^{-is\Delta}F(s)\,ds\|_{L^2_x}\|G\|_{L^1_tL^2_x},
\end{equation}
and \begin{equation}\label{unouno}\|G\|_{L^1_tL^2_x}\lesssim
2^{j/2}\|G\|_{W(L^2,L^{4/3})_tL^2_x}.
\end{equation}
For $a\geq 2$, let now $\tilde{q}=\tilde{q}(a)$ be defined by
\begin{equation}\label{qa}
\frac{2}{\tilde{q}(a)}+\frac{d}{a}=\frac{d}{2}.
\end{equation}
The non-endpoint case of \eqref{dh}, written for $\tilde{r}=a$ and
the $\tilde{q}$ above, gives
\begin{eqnarray*}\sup_t\!\|\int_{t-2^{j+1}<s\leq t-2^j}\!\!\!
e^{-is\Delta}F(s)\,ds\|_{L^2_x}\!\!&=&\!\!\sup_t\|\int_{\R}
e^{-is\Delta}(T_{-t}\chi_{[-2^{j+1},-2^{j}]})(s)F(s)\,ds\|_{L^2_x}\\
&\lesssim& \|F\|_{W(L^{(\tilde{q}/2)'},L^{\tilde{q}})W(\cF
L^a,L^{a'})},
\end{eqnarray*}
for every $2\leq a<r$. In what follows we apply H\"older's
inequality with the triple of indices $2,p,(\tilde{q}/2)^\prime$, so that $\frac{1}{2}+\frac{1}{p}=\frac{1}{(\tilde{q}/2)^\prime}$. This gives, for $g\in\mathcal{C}^\infty_0(\mathbb{R})$,
\begin{eqnarray*}\|FT_zg\|_{L^{(\tilde{q}/2)'}_t(W(\cF
L^a,L^{a'})_x)}&\lesssim&\|FT_zg\|_{L^{(\tilde{q}/2)'}_t(W(\cF
L^a,L^{a'})_x)}\\
&\lesssim&\|FT_zg\|_{L^{2}_t(W(\cF
L^a,L^{a'})_x)}(2^j)^{1/p}\\
&=&\|FT_zg\|_{L^{2}_t(W(\cF L^a,L^{a'})_x)}2^{-j(d/2-d/a-1/2)},
\end{eqnarray*}
where we used
$$\frac{1}{p}=\frac{1}{(\tilde{q}/2)^\prime}-\frac{1}{2}=1-\frac{2}{\tilde{q}}-\frac{1}{2}=\frac12-\frac
d2+\frac da.
$$
Since the support of $F$ with respect to the time is contained in
an interval of duration $\asymp 2^j\leq1 $, the support of the
function
 $\mathbb{R} \ni z\mapsto \|F\,T_z g\|_{L^{2}_t(W(\cF L^a,L^{a'})_x)}$ is contained in an interval
   of duration $\asymp 1$. This allows us to apply H\"older's inequality with respect to the global component,
   too, and we end up with
\begin{align*}\|F\|_{W(L^{(\tilde{q}/2)'},L^{\tilde{q}})_tW(\cF
L^a,L^{a'})_x}&\lesssim&
2^{-j(d/2-d/a-1/2)}\|\|FT_zg\|_{L^{2}_t(W(\cF
L^a,L^{a'})_x)}\|_{L^{\tilde{q}}}\\
&\lesssim & 2^{-j(d/2-d/a-1/2)}\|\|FT_zg\|_{L^{2}_t(W(\cF
L^a,L^{a'})_x)}\|_{L^{4/3}}.
\end{align*}
This estimate, together with \eqref{unozero} and \eqref{unouno}, yields the
estimate \eqref{(ii)}.\par\noindent {\it Case (iii).} Use the same
arguments as in case (ii).
\end{proof}

Since $\Fur L^{r,2}\hookrightarrow \Fur L^r$ and in view of \eqref{primolemma} and \eqref{dec},  in order to prove
\eqref{aim} it suffices to prove
\[
\sum_{j\leq0}|T_j(F,G)|\lesssim \|F\|_{W(L^2,L^{4/3})_t W(\Fur
L^{{r},2},L^{{r}^\prime})_x} \|G\|_{W(L^2,L^{4/3})_t W(\Fur
L^{{r},2},L^{{r}^\prime})_x}.
\]
Now, this can be achieved from \eqref{lammasecondo} by the same
interpolation arguments as in \cite[Par. 6]{keel}. Precisely, we
take $a_0,a_1,b_0,b_1$ such that $(1/r,1/r)$ is inside a small
triangle with vertices $(1/a_0,1/b_0)$, $(1/a_1,1/b_0)$ and
$(1/a_0,1/b_1)$ (see Figure 2), so that
\[
\beta(a_0,b_1)=\beta(a_1,b_0)\not=\beta(a_0,b_0).
\]
\vspace{-0.2cm}
  \begin{center}
           \includegraphics{figstr2.1}
            \\
           $ $
\end{center}
           \begin{center}{ Figure 2.  }
           \end{center}
  \vspace{1.2cm}
Then, we apply Lemma 6.1 of \cite{keel} with $T=\{T_j\}$ (upon
setting $T_j=0$ for $j>0$), $p=q=2$, $r=1$, $C_0=l_\infty^{\beta(a_0,b_0)}$, $C_1=l_\infty^{\beta(a_0,b_1)}$ and,  for $k=0,1$, we take
$$A_k=W(L^2(W(\Fur
L^{a_k},L^{{a_k}^\prime})_x),L^{4/3}),\  B_k=W(L^2(W(\Fur
L^{b_k},L^{{b_k}^\prime})_x),L^{4/3}). $$
Here we choose  $\theta_0$, $\theta_1$, so that
$$1/r=(1-\theta_0)/a_0+\theta_0/a_1, \quad 1/r=(1-\theta_1)/b_0+\theta_1/b_1,$$
(recall, $r=2d/(d-1)$). This gives at once the
desired result, since $\beta(r,r)=0$ and
$$ W(L^2(W(\Fur L^{r,2},L^{{r}^\prime})_x),L^{4/3})\subset
(A_0,A_1)_{\theta_0,2}\cap (B_0,B_1)_{\theta_1,2},$$ as one sees by
applying (in order) Proposition \ref{inter9}, Theorem 1.18.4 of
\cite{triebel}, page 129 (with $p=p_0=p_1=2$), and again
Proposition \ref{inter9} to the Wiener spaces with respect to
$x$.  This concludes the proof of \eqref{hombis} and
\eqref{dhbis}.\par The corresponding retarded estimates can be obtained as follows. The case $(\tilde{q},\tilde{r})=(q,r)=P$ is exactly \eqref{aim}. The case $(\tilde{q},\tilde{r})=P$, $(q,r)\not=P$,  can be obtained by a repeated use of H\"older's inequality to interpolate from the case $(\tilde{q},\tilde{r})=(q,r)=P$, and the case $(\tilde{q},\tilde{r})=P$, $(q,r)=(\infty,2)$ (that is clear from \eqref{dhbis}). Precisely, we want to prove
\[
\|\int_{s<t} e^{i(t-s)\Delta} F(s)\,ds\|_{W(L^{q/2},L^q)_tW(\Fur L^{r^\prime},L^r)_x}\lesssim\|F\|_{W(L^2,L^{4/3})_t W(\Fur L^{\tilde{r},2},L^{\tilde{r}^\prime})},
\]
for $\tilde{r}=2d/(d-1)$, $2/q+d/r=d/2$, $q>4$.\\
We know that such an estimate holds for $(q,r)=(\infty,2)$, as well as
\[
\|\int_{s<t} e^{i(t-s)\Delta} F(s)\,ds\|_{W(L^{2},L^4)_tW(\Fur L^{\tilde{r}^\prime,2},L^{\tilde{r}})_x}\lesssim\|F\|_{W(L^2,L^{4/3})_t W(\Fur L^{\tilde{r},2},L^{\tilde{r}^\prime})}.
\]
Hence, upon setting $I=\int_{s<t} e^{i(t-s)\Delta} F(s)\,ds$, it suffices to prove that
\begin{equation}\label{190}
\|I\|_{W(L^{q/2},L^q)_tW(\Fur L^{r^\prime},L^r)_x}\leq \|I\|^{1-\theta}_{W(L^{2},L^4)_tW(\Fur L^{\tilde{r}^\prime,2},L^{\tilde{r}})_x}\|I\|^\theta_{L^\infty_tL^2_x},
\end{equation}
with $1/r=(1-\theta)/\tilde{r}+\theta/2$, $1/q=(1-\theta)/4+\theta/\infty$, $0<\theta<1$.
To this end, we start with the inequality
\begin{equation}\label{192}
\|u\|_{L^{r'}}\lesssim\|u\|^{1-\theta}_{L^{\tilde{r}^\prime,2}}\|u\|^\theta_{L^{2}},
\end{equation}
(which follows, e.g., from Theorem 1.3.3(g) of \cite{triebel}, page 25, since $(L^{\tilde{r}',2},L^2)_{\theta,1}=L^{r',1}\hookrightarrow L^{r'}$).\par
We apply \eqref{192} with $u=\widehat{I(t)T_zg}$, where $g\in\mathcal{C}^\infty_0(\mathbb{R}^d_x)$ is a non-zero window. We obtain
\begin{align}
\|I(t)\|_{W(\Fur L^{r^\prime},L^r)_x}&=\|\|\widehat{I(t)T_zg}\|_{L^{r'}_x}\|_{L^r_z}\nonumber\\
&\leq\|\|\widehat{I(t)T_zg}\|^{1-\theta}_{L^{{\tilde{r}}',2}_x}\|\widehat{I(t)T_zg}\|^\theta_{L^2_x}\|_{L^r_z}\nonumber\\
&\leq \|I(t)\|^{1-\theta}_{W(\Fur L^{\tilde{r}^\prime,2},L^{\tilde{r}})_x}\|I(t)\|^\theta_{L^2_x},\nonumber
\end{align}
where in the last inequality we applied H\"older's inequality with respect to $z$.
Hence, given any non-zero window $h\in \mathcal{C}^\infty_0(\mathbb{R}_t)$, it follows that
\begin{equation}\label{191}
\|I(t)\|_{W(\Fur L^{r^\prime},L^r)_x}|(T_sh)(t)|\leq\left(\|I(t)\|_{W(\Fur L^{\tilde{r}^\prime,2},L^{\tilde{r}})_x}|(T_sh)(t)|\right)^{1-\theta}\left(\|I(t)\|_{L^2_x}|(T_sh)(t)|\right)^\theta.
\end{equation}
We then take the $L^q_sL^{q/2}_t$ norm of the expression in the left hand side of \eqref{191}
and we apply again H\"older's inequality. This gives the desired result \eqref{190}.\par
Eventually one can prove the retarded estimate in the case $(q,r)=P$, $(\tilde{q},\tilde{r})\not=P$ showing, by the arguments above, that the dual inequality holds for the adjoint operator $
G\mapsto \int _{t>s} e^{-i(t-s)\Delta} G(t)\,dt.$\par
This concludes the proof of Theorem \ref{endpoint}.
\section{An application to Schr\"odinger equations with time-dependent potentials}\label{section6}

Consider the Cauchy problem
\begin{equation}\label{cpP}
\begin{cases}
i\partial_t u+\Delta u=V(t,x)u,\quad t\in [0,T]=I_T, \,\,x\in\rd,\ (d\geq2),\\
u(0,x)=u_0(x),
\end{cases}
\end{equation}
with a potential,
\begin{equation}\label{pot2}V\in L^\alpha(I_T; L^p_x),\quad\frac1{\a}+\frac{d}{2p}\leq1,\ 1\leq\alpha<\infty,\ \frac{d}{2}< p\leq\infty.
\end{equation}
It is proved in \cite[Theorem 1.1, Remark 1.3]{Dancona05} (see also \cite[Theorem 1.1]{Yajima87})  that under the assumpion \eqref{pot2} the Cauchy problem \eqref{cpP}
is well-posed in $L^2$ and admits a unique solution $u\in \mathcal{C}(I_T,\lrd)\cap L^{q}(I_T, L^r)$, for all Schr\"odinger admissible pairs $(q,r)$.\par
We now illustrate an application of our estimates, by deducing a similar result with the solution being controlled in terms of Wiener amalgam norms. We consider the subclass of potentials
\begin{equation}\label{pot3}
V\in L^\alpha(I_T; L^p_x),\quad\frac1{\a}+\frac{d}{p}\leq1,\ 1\leq\alpha<\infty,\ {d}< p\leq\infty.
\end{equation}
Here is our result.
\begin{theorem}
Assume \eqref{pot3}. Then the Cauchy problem \eqref{cpP} has a unique solution $u\in \mathcal{C}(I_T,\lrd)\cap L^{q/2}(I_T, W(\cF L^{r^\prime}, L^r))\cap L^{2}(I_T, W(\cF L^{2d/(d+1),2}, L^{2d/(d-1)}))$, for all $(q,r)$ such that $2/q+d/r=d/2$, $q>4,r\geq2$.
\end{theorem}
\begin{proof}
The proof closely follows the one of \cite[Theorem 1.1, Remark 1.3]{Dancona05} (based on the classical Strichartz estimates). However we present at least the main steps of the proof for the convenience of the reader who is not familiar with that result.\par
First of all, since the interval $I_T$ is bounded, by H\"older's inequality it suffices to prove the theorem when $1/\alpha+d/p=1$. Hence we assume this, and we prove the case $2\leq \alpha<\infty$ and $1\leq\alpha<2$ separately.\\
Let $J=[0,\delta]$ be a small time interval and set, for $q\geq2$, $q\not=4$, $r\geq 1$,
\[
Z_{q/2,r}=L^{q/2}(J;W(\cF L^{r^\prime}, L^r)_x),
\]
\[
Z_{2,2d/(d-1)}=L^{2}(J, W(\cF L^{2d/(d+1),2}, L^{2d/(d-1)})_x),
\]
and $Z=\mathcal{C}(J;L^2)\cap Z_{2,2d/(d-1)}$, with the norm $\|v\|_Z=\max\{\|v\|_{\mathcal{C}(J;L^2)},\|v\|_{Z_{2,2d/(d-1)}}\}$. Notice that, by the arguments at the end of Section \ref{section5}, we have $Z\subset Z_{q/2,r}$ for all $(q,r)$ as in the statement of the theorem. Consider now the integral formulation of the Cauchy problem, namely $u=\Phi(u)$, where
\[
\Phi(v)=e^{it\Delta}u_0+\int_0^t e^{i(t-s)\Delta}V(s)v(s)\,ds.
\]
From Theorems \ref{prima} and \ref{endpoint} it is easy to see that the following estimates hold:
\begin{equation}\label{nan}
\|\Phi(v)\|_{Z_{q/2,r}}\leq C_0 \|u_0\|_{L^2}+C_0\|Vv\|_{Z_{(\tilde{q}/2)',\tilde{r}'}},
\end{equation}
for all $(q,r)$ such that $2/q+d/r=d/2$, $q\geq 4,r\geq2$, and similarly for $(\tilde{q},\tilde{r})$. Since $1/\alpha+d/p=1$ and $\alpha\geq 2$, among such pairs there is a pair $(\tilde{q},\tilde{r})$ such that $1/(\tilde{q}/2)=1/2-1/\alpha$ and $1/\tilde{r}=(d+1)/(2d)-1/p$. We choose such a pair in \eqref{nan} and, after using the inclusion $L^{(\tilde{q}/2)'}(J;L^{\tilde{r}'})\subset Z_{(\tilde{q}/2)',\tilde{r}'}$ and $Z_{2,2d/(d-1)}\subset L^2(J;L^{2d/(d-1)})$, we apply H\"older's inequality. We obtain
\begin{equation}\label{nan2}
\|\Phi(v)\|_{Z_{q/2,r}}\leq C_0 \|u_0\|_{L^2}+C_0\|V\|_{L^\alpha(J; L^p)}\|v\|_{Z_{2,2d/(d-1)}}.
\end{equation}
By taking $(q,r)=(\infty,2)$ and $(q,r)=(4,2d/(d-1))$ one deduces that $\Phi:Z\to Z$ (the fact that $\Phi(u)$ is {\it continuous} in $t$ when valued in $L^2_x$ follows from a classical limiting argument). Also, since $\alpha<\infty$, if $J$ is small enough, $C_0\|V\|_{L^\alpha_t L^p_x}<1/2$, and $\Phi$ is a contraction. This gives a unique solution in $J$. By iterating this argument a finite number of times one obtains a solution in $[0,T]$.\par
In the case $1\leq \alpha<2$, one starts instead from the estimate
\[
\|\Phi(v)\|_{Z_{q/2,r}}\leq C_0\|u_0\|_{L^2}+C_0\|Vv\|_{Z_{\alpha,2p/(p+2)}}
\]
which follows again from Theorems \ref{prima} and \ref{endpoint} with $((\tilde{q}/2)',\tilde{r}')=(\alpha,2p/(p+2))$. In fact, from $1\leq\alpha<2$ and $1/\alpha+d/p=1$, it follows $\tilde{q}>4$ and $2/\tilde{q}+d/\tilde{r}=d/2$.
Then, again by H\"older's inequality, one obtains
\[
\|\Phi(v)\|_{Z_{q/2,r}}\leq C_0\|u_0\|_{L^2}+C_0\|V\|_{L^\alpha(J;L^p)}\|v\|_{Z_{\infty,2}}.
\]
From here, proceeding  as above yields  the desired result.

\end{proof}

\vskip0.5truecm
\textit{Acknowledgments}. The authors thank Luigi Rodino  for very helpful  discussions on this topic.


\begin{thebibliography}{10}
\bibitem{baoxiang} W. Baoxiang, Z. Lifeng and G. Boling. Isometric decomposition operators, function spaces $E_{p,q}^\lambda$ and applications to nonlinear evolution equations. {\it J. Funct. Anal.}, 233(1):1--39, 2006.
\bibitem{benyi} A. B\'enyi, K. Gr\"ochenig, K.A. Okoudjou and L.G. Rogers. Unimodular Fourier multipliers for modulation spaces. Preprint, 2006.

\bibitem{Dancona05} P. D'Ancona, V. Pierfelice  and N. Visciglia.
             Some remarks on the {S}chr\"odinger equation with a potential
              in {$L\sp r\sb tL\sp s\sb x$}. {\it Math. Ann.}, 333(2):271--290, 2005.
\bibitem{feichtinger80}
H.~G. Feichtinger.
\newblock Banach convolution algebras of {W}iener's type,
\newblock In {\em Proc. Conf. ``Function, Series, Operators",
 Budapest August 1980}, Colloq. Math. Soc. J\'anos Bolyai, 35,  509--524, North-Holland, Amsterdam,
 1983.
\bibitem{feichtinger83}
H.~G. Feichtinger.
\newblock Banach spaces of distributions of {W}iener's type and interpolation.
\newblock In {\em Proc. Conf. Functional Analysis and Approximation,
 Oberwolfach August 1980},  Internat. Ser. Numer. Math., 69:153--165. Birkh\"auser, Boston, 1981.
\bibitem{feichtinger90}
H.~G. Feichtinger.
\newblock Generalized amalgams, with applications to {F}ourier transform.
\newblock {\em Canad. J. Math.}, 42(3):395--409, 1990.
\bibitem{fournier-stewart85}
J.~J.~F. Fournier and J.~Stewart.
\newblock Amalgams of ${L}\sp p$ and $l\sp q$.
\newblock {\em Bull. Amer. Math. Soc. (N.S.)}, 13(1):1--21, 1985.
\bibitem{Fei98}
H.~G. Feichtinger and G.~Zimmermann.
\newblock {A} {B}anach space of test functions for {G}abor analysis.
\newblock In H.~G. Feichtinger and T. Strohmer, editors, {\em Gabor Analysis and Algorithms. Theory and Applications}, Applied and Numerical Harmonic Analysis, 123--170. Birkh\"auser, Boston, 1998.
\bibitem{folland89}
G.~B. Folland.
\newblock {\em Harmonic Analysis in Phase Space}.
\newblock Princeton Univ. Press, Princeton, NJ, 1989.
\bibitem{GinibreVelo92}
J.~Ginibre and G.~Velo.
\newblock Smoothing properties and retarded estimates for some dispersive
  evolution equations.
\newblock {\em Comm. Math. Phys.}, 144(1):163--188, 1992.
\bibitem{Heil03}
C. Heil.
\newblock {A}n introduction to weighted {W}iener amalgams.
\newblock In M. Krishna, R. Radha and S. Thangavelu, editors, {\em
Wavelets and their Applications}, 183--216. Allied Publishers Private Limited, 2003.
\bibitem{Kato70}
T.~Kato.
\newblock Linear evolution equations of ``hyperbolic'' type.
\newblock {\em J. Fac. Sci. Univ. Tokyo Sect. I}, 17:241--258, 1970.
\bibitem{keel}
M. Keel and T. Tao.
\newblock { Endpoint Strichartz estimates}.
\newblock {\it Amer. J. Math.}, 120:955--980, 1998.
\bibitem{stein}
E. M. Stein.
\newblock {\it Singular integrals and differentiability properties of functions}.
\newblock Princeton University Press, Priceton, 1970.
\bibitem{stein93}
E. M. Stein.
\newblock {\it Harmonic analysis}.
\newblock Princeton University Press, Priceton,1993.
\bibitem{steinweiss}
E. M. Stein and G. Weiss.
\newblock {\it Introduction to Fourier Analysis on Euclidean spaces}.
\newblock Princeton university Press, 1971.
\bibitem{tao2}
T. Tao.
\newblock {Spherically averaged endpoint  Strichartz estimates for the two-dimensional Schr\"odinger equation}. {\it Comm.\ Partial differential Equations}, 25:1471--1485, 2000.
\bibitem{tao}
T. Tao.
\newblock {\it Nonlinear Dispersive Equations: Local and Global Analysis}.
\newblock CBMS Regional Conference Series in Mathematics, Amer. Math. Soc., 2006.
\bibitem{triebel}
H. Triebel.
\newblock {\it Interpolation theory, function spaces, differential operators}.
\newblock North-Holland, 1978.
\bibitem{Yajima87}
K.~Yajima.
\newblock Existence of solutions for {S}chr\"odinger evolution equations.
\newblock {\em Comm. Math. Phys.}, 110(3):415--426, 1987.
\end{thebibliography}
\end{document}